\pdfoutput=1
\RequirePackage{ifpdf}
\ifpdf 
\documentclass[pdftex]{sigma}
\else
\documentclass{sigma}
\fi

\DeclareMathOperator{\ad}{ad}
\DeclareMathOperator{\psupp}{psupp}
\DeclareMathOperator{\supp}{supp}

\begin{document}

\allowdisplaybreaks

\renewcommand{\thefootnote}{$\star$}

\newcommand{\arXivNumber}{1405.7122}

\renewcommand{\PaperNumber}{115}

\FirstPageHeading

\ShortArticleName{The Freiheitssatz for Generic Poisson Algebras}

\ArticleName{The Freiheitssatz for Generic Poisson Algebras\footnote{This paper is a~contribution to the Special Issue
on Poisson Geometry in Mathematics and Physics.
The full collection is available at \href{http://www.emis.de/journals/SIGMA/Poisson2014.html}
{http://www.emis.de/journals/SIGMA/Poisson2014.html}}}

\Author{Pavel S.~KOLESNIKOV~$^\dag$, Leonid G.~MAKAR-LIMANOV~$^\ddag$ and Ivan P.~SHESTAKOV~$^{\S\dag}$}

\AuthorNameForHeading{P.S.~Kolesnikov, L.G.~Makar-Limanov and I.P.~Shestakov}

\Address{$^\dag$~Sobolev Institute of Mathematics, 630090, Novosibirsk, Russia}
\EmailD{\href{mailto:pavelsk@math.nsc.ru}{pavelsk@math.nsc.ru}}

\Address{$^\ddag$~Department of Mathematics, Wayne State University, Detroit, MI 48202, USA}
\EmailD{\href{mailto:lml@math.wayne.edu}{lml@math.wayne.edu}}

\Address{$^\S$~Instituto de Matematica e Estat\'\i stica, Universidade de S\~ao Paulo, 05508-090 S\~ao Paulo, Brasil}
\EmailD{\href{mailto:shestak@ime.usp.br}{shestak@ime.usp.br}}

\ArticleDates{Received July 29, 2014, in f\/inal form December 22, 2014; Published online December 29, 2014}

\Abstract{We prove the Freiheitssatz for the variety of generic Poisson algebras.}

\Keywords{Freiheitssatz; Poisson algebra; generic Poisson algebra; algebraically closed al\-gebra; polynomial identity;
dif\/ferential algebra}

\Classification{17A30; 17B63; 17A50}

\rightline{\it Dedicated to the memory of Myung Hyo Chul (1937--2010)}

\renewcommand{\thefootnote}{\arabic{footnote}}
\setcounter{footnote}{0}

\section{Introduction}

One of the classical achievements of the combinatorial group theory is the decidability of the word problem in
a~f\/initely generated group with one def\/ining relation~\cite{Magnus}.
This result was a~corollary of a~fundamental statement called Freiheitssatz: every equation over a~free group is
solvable in some extension.
For solvable and nilpotent groups, this complex of problems was studied in~\cite{Rom72}.

In the context of Lie algebras, similar statements were proved~\cite{Sh62}.
For associative algebras, the problem turns to be surprisingly dif\/f\/icult: over a~f\/ield of characteristic zero, the
Freiheitssatz was proved in~\cite{LML85}, but the question about decidability of the word problem for an associative
algebra with one def\/ining relation remains open.
One of the reasons for this is that the variety of associative algebras, contrary to those of groups or Lie algebras, is
not a~Schreier one, that is, a~subalgebra of a~free associative algebra is not necessary free.
And as a~matter of fact, the free algebras of Schreier varieties are usually more easy to deal with.

In~\cite{KMLU08}, the Freiheitssatz was proved for right-symmetric (pre-Lie) algebras, and in~\cite{MLU10}~-- for
Poisson algebras over a~f\/ield of zero characteristic.

There is plenty of varieties for which the Freiheitssatz is not true, e.g., so is the variety of Poisson algebras over
a~f\/ield of positive characteristic.
One may f\/ind more examples of this kind in~\cite{MikSh2013}, e.g., for Leibniz algebras the Freiheitssatz does not hold
(as well as for every variety of di-algebras in the sense of~\cite{Kol2008}).

Throughout the paper $\Bbbk $ denotes a~f\/ield of characteristic zero.

A generic Poisson algebra ($\operatorname{GP}$-algebra) is a~linear space with two operations and one constant:
\begin{itemize}\itemsep=0pt
\item associative and commutative product $x\cdot y = xy$;
\item anti-commutative bracket $\{x,y \}$;
\item multiplicative identity~$1$, $x\cdot 1 = 1\cdot x = x$,
\end{itemize}
satisfying the Leibniz identity
\begin{gather*}
\{x,yz\} = \{x,y\}z + \{x,z\}y.
\end{gather*}
These algebras were introduced in~\cite{Shest2000} in the study of speciality and deformations of Malcev--Poisson
algebras.

Let $\operatorname{AC}(X)$ be the free anti-commutative algebra ($\operatorname{AC}$-algebra) generated by a~set~$X$ with respect to operation denoted
by $\{\cdot,\cdot\}$, and let $\operatorname{GP}(X)$ be the free $\operatorname{GP}$-algebra with a~set of generators~$X$.
As a~linear space, $\operatorname{GP}(X)$ is isomorphic to the symmetric algebra $S(\operatorname{AC}(X))$~\cite{Shest1993}.

There are two general ways to prove the Freiheitssatz for a~particular variety: one may either analyze the structure of
a~one-generated ideal in the free algebra (this often leads to a~solution of the word problem for algebras with one
def\/ining relation) or try to construct an example of an algebra in the considered variety which is ``suf\/f\/iciently rich''
but ``relatively free''.
The last approach was used in~\cite{LML85} and~\cite{MLU10}.
In this paper, we formalize these notions in a~precise form (Section~\ref{sec1}) for an (almost) arbitrary variety.
This may be helpful for study Freiheitssatz in other varieties.
By ``suf\/f\/iciently rich'' algebraic system we mean an algebra which is algebraically closed in certain sense
(conditionally closed).
Conditionally closed systems may not be algebraically closed in the strong sense of~\cite{Bokut66} or even existentially
closed (see, e.g.,~\cite{Higman}), but they turn to be extremely helpful for studying one-generated ideals of free
algebras.
Our construction of a~conditionally closed generic Poisson algebra is completely dif\/ferent from the one in~\cite{MLU10}.
We apply general results on equations in dif\/ferential algebras and then twist the standard functor relating dif\/ferential
algebras and Poisson algebras to obtain a~generic Poisson algebra which is ``relatively free''.
For an algebra~$A$, the property to be ``relatively free'' means the existence of a~suf\/f\/iciently large lower bound for
the degree of a~polynomial identity holding on~$A$.
For example, such an algebra may have no nontrivial polynomial identities at all.

The proof that the system obtained is indeed ``relatively free'' is inspired by~\cite{FarkasI}, but the reali\-za\-tion of
the scheme has one distinguishing feature.
As in~\cite{FarkasI}, we show that if a~generic Poisson algebra has a~polynomial identity then it has an identity of
a~special form.
The proof is split into three main steps: (1)~classif\/ication of anti-commutative polynomials that are derivations in one
variable; (2)~classif\/ication of anti-commutative polynomials that are derivations in all variables (Jacobian
polynomials); (3)~classif\/ication of Jacobian generic Poisson polynomials.
The f\/irst and third steps do not essentially dif\/fer from~\cite{FarkasI}, but the second step (which is almost trivial
for Lie algebras) is highly nontrivial in the generic (anti-commutative) case.
In particular, there are no (poly-linear) Lie polynomials that are derivations with respect to more than two variables.
In contrast, one may invent an anti-commutative polynomial which is a~derivation with respect to any number of
variables.
However, it turns out that just a~limited number of anti-commutative polynomials may be Jacobian.
We describe all such polynomials in Section~\ref{sec2}.

\section{Conditionally closed algebras and the Freiheitssatz}\label{sec1}

Suppose $\mathfrak M$ is a~variety of algebras over a~f\/ield~$\Bbbk $.
Denote by $\mathfrak M(X)$ the free algebra in $\mathfrak M$ generated by a~set~$X$.
For $A,B\in \mathfrak M$, the notation $A*_{\mathfrak M} B$ stands for the free product of~$A$ and~$B$ in $\mathfrak M$.

If $A\in \mathfrak M$ then every $\Psi \in A*_{\mathfrak M} \mathfrak M(x)$ may be interpreted as an~$A$-valued
function on~$A$.
Moreover, for every extension $\bar A$ of~$A$, $\bar A\in \mathfrak M$, $\Psi(x)$ is an $\bar A$-valued function
on~$\bar A$.
An equation of the form $\Psi(x)=0$ is {\em solvable over}~$A$ if there exists an extension $\bar A\in \mathfrak M$
of~$A$ such that the equation has a~solution in $\bar A$.
If such a~solution can be found in~$A$ itself then $\Psi(x)=0$ is said to be {\em solvable in}~$A$.

Recall the common def\/inition (see, e.g.,~\cite{Higman, Scott}): an algebra~$A$ is {\em $($existentially$)$ algebraically
closed\/} if every system of equations which is solvable over~$A$ is solvable in~$A$.
Let us restrict this def\/inition to a~particular case of one equation: we will say $A\in \mathfrak M$ to be existentially
closed in~$\mathfrak M$ if every equation $\Psi(x)=0$, where $\Psi(x)\in A*_{\mathfrak M} \mathfrak M(x)$, which is solvable in an
appropriate extension~$\bar A$ of~$A$, $\bar A\in \mathfrak M$, has a~solution in~$A$.
This def\/inition is important for model theory, and it can be an ef\/f\/icient tool for studying algebras provided the
principal question on the solvability of a~particular equation is solved.

A stronger property (see~\cite{Bokut66}) can be stated as follows: an algebra $A\in \mathfrak M$ is called {\em
algebraically closed in $\mathfrak M$} if for every $\Psi \in A*_{\mathfrak M} {\mathfrak M}(x)$, $\Psi\notin A$, the
equation $\Psi(x)=0$ is solvable in~$A$.
We are going to propose an intermediate def\/inition which is suf\/f\/icient for our purpose.

\begin{definition}
An algebra $A\in \mathfrak M$ is called {\em conditionally closed in $\mathfrak M$} if for every $\Psi(x)\in
A*_{\mathfrak M} {\mathfrak M}(x)$ which is not a~constant function on~$A$ the equation $\Psi(x)=0$ is solvable in~$A$.
\end{definition}

Every algebraically closed in $\mathfrak M$ algebra is conditionally closed in $\mathfrak M$.
However, there is plenty of conditionally closed systems that are not algebraically closed in $\mathfrak M$.
For example, an algebraically closed f\/ield is conditionally closed but not algebraically closed in the variety of all
associative algebras.
Similarly, such a~f\/ield may be considered as a~Poisson algebra with respect to trivial bracket, and the Poisson algebra
obtained is conditionally closed but not algebraically closed in the variety of all Poisson algebras.

It is also interesting to compare conditionally closed and existentially closed algebras.
Neither of these notions is a~formal generalization of another.
For example, let $\mathfrak M=\operatorname{As}$, the variety of associative algebras, and let $A\in \operatorname{As}$ be the algebraic closure of
the f\/ield $\Bbbk(t)$.
As an algebraically closed f\/ield,~$A$ is conditionally closed in $\operatorname{As}$, but the equation $[t,x]=1$ has no solution in~$A$
although it is solvable in an appropriate extension (e.g., in the Makar-Limanov's skew f\/ield~\cite{LML85}).
On the other hand, for the same variety $\operatorname{As}$, the existential algebraic closure (see, e.g.,~\cite[Chapter~III]{Cherlin})
of quaternions $\mathbb H$ is not conditionally closed: equation $ix-xi=1$ has no solution in any extension of $\mathbb
H$ in~$\operatorname{As}$.

Suppose $\mathfrak M_1$ and $\mathfrak M_2$ are two varieties of algebras over a~f\/ield $\Bbbk $, and let $\omega:
\mathfrak M_1 \to \mathfrak M_2$ be a~functor which acts as follows: given $A\in \mathfrak M_1$, $A^{(\omega)}\in
\mathfrak M_2$ is the same linear space equipped with new operations expressed in terms of initial operations.
For example, one may consider the classical functor from the variety of associative algebras into the variety of Lie
algebras def\/ined by $[x,y] = xy-yx$.

Another important example comes from the following settings.
Let $\mathfrak M_1 = \operatorname{Dif}_{2n}$ be the variety of commutative associative algebras with $2n$ pairwise commuting
derivations $\partial_i$, $\partial'_i$, $i=1,\dots, n$.
Then, given $A\in \operatorname{Dif}_{2n}$, the same space equipped with new binary operation
\begin{gather}
\label{eq:PoisBracket}
\{a,b\} = \sum\limits_{i=1}^n \partial_i(a)\partial'_i(b) - \partial_i(b)\partial'_i(a),
\qquad
a,b\in A,
\end{gather}
is known to be a~Poisson algebra denoted by $A^{(\partial)}$.
If we allow the derivations $\partial_i$, $\partial'_i$ to be non-commuting then~\eqref{eq:PoisBracket} def\/ines on
a~commutative algebra~$A$ a~structure of a~$\operatorname{GP}$-algebra.

In general, $\omega $ may be a~functor induced by a~morphism of the governing operads.
Functors of this kind were closely studied in~\cite{MikSh2013}.

\begin{proposition}
\label{prop:CCFunctor}
If an algebra $A\in \mathfrak M_1$ is conditionally closed in $\mathfrak M_1$ then $A^{(\omega)}$ is conditionally
closed in $\mathfrak M_2$.
\end{proposition}

Note that for algebraically closed algebras this statement does not hold.

\begin{proof}
Since $\omega $ is a~functor, the universal property of the free product implies the existence of a~homomorphism
$\varphi: A^{(\omega)}*_{\mathfrak M_2} \mathfrak M_2(x) \to (A*_{\mathfrak M_1} \mathfrak M_1(x))^{(\omega)} $ such
that $f(a)=\varphi(f)(a)$ for all $f=f(x)\in A^{(\omega)}*_{\mathfrak M_2} \mathfrak M_2(x)$, $a\in A^{(\omega)}$.

Therefore,~$f$ is not a~constant function on $A^{(\omega)}$ if and only if $\varphi(f)$ is not a~constant function
on~$A$.
If~$A$ is conditionally closed then there exists $a\in A$ such that $\varphi(f)(a)=0$ and thus $f(a)=0$.
\end{proof}

In some cases, the converse statement is true: if $A^{(\omega)}$ is conditionally closed in $\mathfrak M_2$ then so
is~$A$ in $\mathfrak M_1$.

\begin{remark}
\label{rem:CCAlgebras_const}
Let $\omega: \mathfrak M_1 \to \mathfrak M_2$ be a~functor with the following property: if $A\in \mathfrak M_1$ and
$A^{(\omega)}$ is a~subalgebra of $C\in \mathfrak M_2$ then there exists $B\in \mathfrak M_1$ such that $C=B^{(\omega)}$
and~$A$ is a~subalgebra of~$B$.
Then~$A$ is conditionally closed in $\mathfrak M_1$ provided that $A^{(\omega)}$ is conditionally closed in $\mathfrak
M_2$.
\end{remark}

Indeed, $A^{(\omega)}*_{\mathfrak M_2} \mathfrak M_2(x)$ is a~$\mathfrak M_2$-algebra which contains $A^{(\omega)}$.
Hence, $A^{(\omega)}*_{\mathfrak M_2} \mathfrak M_2(x) = B^{(\omega)}$ for a~$\mathfrak M_1$-algebra~$B$.
Therefore, there exists a~homomorphism of $\mathfrak M_1$-algebras  $\psi: A*_{\mathfrak M_1} \mathfrak M_1(x) $ $\to B $
such that $\psi(a)=a$ for $a\in A$, $\psi(x)=x$.
Hence, $\varphi(\psi(f)) = f $ for all $f\in A*_{\mathfrak M_1} \mathfrak M_1(x)$, where~$\varphi $ is the homomorphism
in the proof of Proposition~\ref{prop:CCFunctor}.

Suppose $A^{(\omega)}$ is conditionally closed in $\mathfrak M_2$.
If $f \in A*_{\mathfrak M_1} \mathfrak M_1(x)$ is not a~constant function on~$A$ then so is $g =\psi(f)$ since
$\varphi(g(a)) = \varphi(g)(a) = f(a)$ for all $a\in A$.
Hence, there exists a~solution of the equation $g(x)=0$ in $A^{(\omega)}$ which is obviously a~solution of $f(x)=0$
in~$A$.

The {\em Freiheitssatz problem} for a~variety $\mathfrak M$ is to determine whether every nontrivial equation over the
free algebra $\mathfrak M(X)$, $X=\{x_1,x_2,\dots \}$, is solvable over $\mathfrak M(X)$.

It is obviously equivalent to the following question about free algebras: is the intersection of the ideal $(f)$
generated by an element $f\in \mathfrak M(X\cup \{x\})$ and the subalgebra $\mathfrak M(X)\subset \mathfrak M(X\cup
\{x\})$ trivial if $f \notin \mathfrak M(X)$ (i.e., depends on~$x$)? If the answer is positive for all such~$f$ then we
say that the Freiheitssatz holds for $\mathfrak M$.

It is easy to note that if the Freiheitssatz holds of a~variety $\mathfrak M$ over algebraically closed f\/ields then it
holds for $\mathfrak M$ over an arbitrary f\/ield of the same characteristic.

\begin{lemma}
\label{lem:CC-Frei}
Suppose $\mathfrak M$ is a~variety of algebras with at least one binary operation $\cdot $ in the signature such that
$\mathfrak M(X)=\mathfrak M(x_1,x_2,\dots)$ has no zero divisors with respect to $\cdot $.
Then, if for every nonzero polynomial $h=h(x_1,\dots, x_n)\in \mathfrak M(X)$ there exists a~conditionally closed
algebra $A\in \mathfrak M$ which does not satisfy the polynomial identity $h(x_1,\dots, x_n) = 0$ then the Freiheitssatz
holds for $\mathfrak M$.
\end{lemma}

\begin{proof}
Suppose $X=\{x_1,x_2,\dots \}$, $x\notin X$, and let $f=f(x, x_1,\dots, x_n)\in \mathfrak M(X\cup\{x\})\setminus{\mathfrak M}(X)$.
Then $f=f_1 + f_0$, where $f_1$ belongs to the ideal generated by~$x$, $f_0\in {\mathfrak M}(X)$.

Assume $g\in (f)\cap {\mathfrak M}(X)$, $g\ne 0$.
Then $h=f_1 \cdot g \ne 0$, hence, there exist a~conditionally closed $A\in \mathfrak M$ such that $h(x,x_1,\dots, x_n)$
is not a~polynomial identity on~$A$.
Therefore, there exist $a,a_1,\dots, a_n\in A$ such that $f_1(a, a_1,\dots, a_n) g(a_1,\dots, a_n)\ne 0 $ in~$A$, so
\mbox{$f_1(a, a_1,\dots, a_n)\ne 0$}.
On the other hand, $f_1(0,a_1,\dots, a_n)=0$.
Therefore, $\Psi(x) = f_1(x, a_1,\dots, a_n)$ is a~non-constant function on~$A$.
Since~$A$ is conditionally closed, there exists $a\in A$ such that $\Psi(a) = f_1(a, a_1,\dots, a_n) = -f_0(a_1,\dots,
a_n)$.
Thus, $f(a,a_1,\dots, a_n)=0$ but $g(a_1,\dots, a_n)\ne 0$ which is impossible if $g\in (f)\triangleleft \mathfrak
M(X\cup\{x\})$.
\end{proof}

There is a~well-known functor $\omega $ from the variety $\operatorname{As}$ of associative algebras to the variety $\operatorname{Jord}$ of Jordan
algebras: every associative algebra~$A$ turns into a~Jordan algebra denoted by~$A^{(+)}$ under new product $x\circ y=xy + yx$.
A~Jordan algebra is said to be {\em special} if it can be embedded into an algebra of the form~$A^{(+)}$, $A\in \operatorname{As}$.
The class of all special Jordan algebras is not a~variety since a~homomorphic image of a~special Jordan algebra may not be special.
However, the class of all homomorphic images of all special Jordan algebras is a~variety denoted by~$\operatorname{SJ}$.
The free algebra~$\operatorname{SJ}(X)$ is obviously the subalgebra of $\operatorname{As}(X)^{(+)}$ generated by~$X$ with respect to Jordan product.

\begin{corollary}
The Freiheitssatz holds for the variety generated by special Jordan algebras over a~field of characteristic zero.
\end{corollary}

\begin{proof}
Let us recall the main features of the construction of the algebraically closed associative noncommutative algebra
from~\cite{LML85}.
Consider the associative commutative algebra~$C$ generated by all rational powers of countably many formal variables
$p_i$, $q_i$, $i\in \mathbb N$:
\begin{gather*}
C = \Bbbk\big[p_i^\lambda, q_i^\mu \,\big|\, \lambda,\mu\in \mathbb Q,\; i=1,2,\dots\big].
\end{gather*}
The algebra~$C$ has two derivations $\partial_p$ and $\partial_q$ def\/ined by inductive relations
\begin{gather*}
\partial_p(p_{n+1}) = p_n^{-1}\partial_p(p_n),
\qquad
\partial_p(q_{n+1}) = p_n^{-1}\partial_p(p_n) + q_n^{-1}\partial_p(q_n),
\\
\partial_p(q_{n+1}) = 0,
\qquad
\partial_q(q_{n+1}) = q_n^{-1}\partial_q(q_n),
\end{gather*}
assuming $\partial_p(p_1)=\partial_q(q_1)=1$, $\partial_p(q_1)=\partial_q(p_1)=0$.

The set of all monomials in~$C$ is an Abelian group which is linearly ordered by means of
\begin{gather*}
p_1 \ll q_1 \ll p_2 \ll q_2 \ll \dots <1,
\end{gather*}
where $a\ll b$ means $a^\lambda <b$ for all $\lambda >0$.
Let us denote this group by~$G$, and let $G_{n,m}$ stand for the subgroup of~$G$ generated by $ p_i^\lambda $, $q_i^\mu
$, $i=1,\dots, n$, $\lambda,\mu \in \frac{1}{m}\mathbb Z$.

A Malcev--Newmann series over~$G$ is a~transf\/inite formal series
\begin{gather*}
a = \sum\limits_{g\in G} a(g) g,
\qquad
a(g)\in \Bbbk,
\end{gather*}
in which the set $\supp(a) = \{g\in G\,|\, a(g)\ne 0\}$ is a~well-ordered subset of~$G$.

Denote by $A_{n,m}$ the set of all Malcev--Newmann series~$a$ over~$G$ with $\supp(a)\subset G_{n,m}$, $n,m\in
\mathbb N$.
The union~$A$ of all $A_{n,m}$ is an associative and commutative subalgebra in the algebra of all Malcev--Newmann series
over~$G$.
Derivations~$\partial_p$ and~$\partial_q$ may be naturally expanded to derivations of~$A$.

Finally, def\/ine a~new binary operation $*$ on the space~$A$:
\begin{gather*}
a*b = \sum\limits_{n=0}^\infty \frac{1}{n!} \partial_q^n(a) \partial_p^n(b).
\end{gather*}
It is proved in~\cite{LML85} that~$A$ with respect to $*$ is an algebraically closed skew f\/ield.
Note that
\begin{gather*}
q_1*p_1 - p_1*q_1 = 1,
\end{gather*}
so~$A$ contains the f\/irst Weyl algebra $W_1$.
Since~$A$ is a~skew f\/ield, it contains the skew f\/ield of fractions~$Q(W_1)$, which contains a~two-generated free
associative algebra~\cite{LML83}.
Thus,~$A$ contains free associative algebra in any f\/inite number of generators $x_1,\dots, x_n$.
The special Jordan algebra~$A^{(+)}$ is conditionally closed by Proposition~\ref{prop:CCFunctor} and contains free
special Jordan algebra $\operatorname{SJ}(x_1,\dots, x_n)$.
Therefore, the variety $\operatorname{SJ}$ satisf\/ies all conditions of Lemma~\ref{lem:CC-Frei}.
\end{proof}

Note that for the entire variety $\operatorname{Jord}$ the Freiheitssatz does not hold~\cite{MikSh2013}.

\section{Jacobian polynomials in free anti-commutative algebras}

In order to prove the Freiheitssatz for a~variety $\mathfrak M$ by means of Lemma~\ref{lem:CC-Frei}, we have to
construct an algebra in $\mathfrak M$ which is conditionally closed  and does not satisfy a~given polynomial
identity.

In this section, we discuss technical questions that are used in subsequent sections for the study of polynomial
identities on generic Poisson algebras.

\subsection[Preliminaries on $\operatorname{AC}(X)$]{Preliminaries on $\boldsymbol{\operatorname{AC}(X)}$}

Let~$X$ be a~set of generators, and let $X^*$ stand for the set of all (nonempty) associative words~$u$ in the alphabet~$X$.
Denote by $X^{**}$ the set of all non-associative words in~$X$.
Given a~word $u\in X^*$, denote by $(u)$ a~non-associative word obtained from~$u$ by some bracketing.
We will also use $[X^*]$ to denote the set of all associative and commutative words in~$X$.
Given $u\in X^*$, $[u]$~stands for the commutative image of~$u$.

Suppose $X^{**}$ is equipped with a~linear order $\preceq $.
A~non-associative word $u\in X^{**}$ is {\em normal} if either $u=x\in X$ or $u=u_1u_2$, where $u_1$ and $u_2$ are
normal and $u_1\prec u_2$.
Obviously, normal words in $X^{**}$ form a~linear basis of the free anti-commutative algebra $\operatorname{AC}(X)$ generated by~$X$
(see~\cite{BokutCL, ShirshovE}).

Let us call the elements of $\operatorname{AC}(X)$ {\em $\operatorname{AC}$-polynomials}.
Given $u\in X^{**}$, def\/ine $\deg u$ to be the length of~$u$.
Thus, we have a~well-def\/ined degree function on $\operatorname{AC}(X)$.

Choose a~generator $x_i\in X=\{x_1,\dots, x_n\}$ and denote by $V_i$ the subspace of $\operatorname{AC}(X)$ spanned by all
nonassociative words linear in~$x_i$.
Fix a~linear order $\preceq $ on $X^{**}$ such that any nonassociative word which contains $x_i$ is greater than any
word without $x_i$ (there exist many linear orders with this property).
With respect to such an order, the unique normal form of a~monomial $w\in V_i$ is
\begin{gather}
\label{eq:Normal-i-Form}
w = \{u_1, \{u_2, \dots \{u_k, x_i\}\dots\}\},
\end{gather}
where $u_j$, $j=1,\dots, k$, are normal words in the alphabet $X\setminus \{x_i\}$.
The number~$k$ is called {\em $x_i$-height} \cite{FarkasI} of~$w$, let us denote it by $\operatorname{ht}(w, x_i)$.

Let $V_0 = \bigcap\limits_{i=1}^n V_i$ be the space of polylinear $\operatorname{AC}$-polynomials.
It is easy to compute $x_i$-height of a~nonassociative word $w\in V_0$ just by the number of brackets in~$w$ to the left
of~$x_i$, assuming~$\{$ is counted as~1 and $\}$ as~$-1$.
For example, the $x_4$-height of $\{\{x_1, \{\{x_2, x_3\},x_4\} \},\{x_5,x_6\}\}$ is equal to~3.

Denote by $M(\operatorname{AC}(X))$ the algebra of left multiplications on $\operatorname{AC}(X)$, i.e., the subalgebra of $\mathrm{End}_\Bbbk
\operatorname{AC}(X)$ generated~by
\begin{gather*}
\ad g:  \ f\mapsto \{g,f\},
\qquad
f,g\in \operatorname{AC}(X).
\end{gather*}
Since the variety $\operatorname{AC}$ of anti-commutative algebras is a~Schreier one, $M(\operatorname{AC}(X))$ is a~free associative algebra
(see~\cite{Umirbaev1994}).
Let~$U$ stand for the set of all normal words in $X^{**}$.
It is easy to see that $M(\operatorname{AC}(X)) \simeq  \operatorname{As}(U)$ provided that we identify $\ad u$ with $u\in U$.

Denote by $U_i$ the set of normal words in the alphabet $X\setminus \{x_i\} $.
Then $V_i$ is a~1-generated free left module over $\operatorname{As}(U_i)$: every word of the form~\eqref{eq:Normal-i-Form} may be
uniquely presented as
\begin{gather*}
w = W(x_i),
\qquad
W=\ad u_1 \ad u_2 \cdots \ad u_k,
\end{gather*}
where $u_1,\dots, u_k \in U_i$.

Denote by $*$ the involution of $\operatorname{As}(U_i)$ given by $(u_1\cdots u_k)^* = (-1)^{k}u_k\cdots u_1$, $u_j\in U_i$.

\begin{definition}
\label{defn:x_i-flip}
A~linear transformation of $V_i$ def\/ined by the rule
\begin{gather*}
F_i: \ W(x_i) \mapsto -W^*(x_i)
\end{gather*}
is called an {\em $x_i$-flip}.
Obviously, $(F_i)^{-1}=F_i$.
\end{definition}

The set of all f\/lips $\{F_1,\dots, F_n\}$ acts on the space $V_0$ and thus generates a~group $\mathcal F\subseteq
\mathrm{GL}(V_0)$.
Given a~normal word $u\in V_0$, the orbit $\mathcal Fu$ consists of $\operatorname{AC}$-monomials (polynomials of the form~$\varepsilon v$, $v$~is a~nonassociative word, $\varepsilon =\pm 1$).

\begin{lemma}
\label{lem:FlipsGenerateS_n}
Let $w = \{x_1,\{x_2, \dots \{x_{n-1},x_n\}\dots \}\}\in V_0$.
Then
\begin{gather*}
(-1)^\sigma \{x_{1\sigma},\{x_{2\sigma}, \dots \{x_{(n-1)\sigma},x_{n\sigma}\}\dots \}\}\in \mathcal Fw
\end{gather*}
for every $\sigma \in S_n$ $($here $(-1)^\sigma $ stands for the parity of a~permutation~$\sigma)$.
\end{lemma}

\begin{proof}
In order to apply $F_i$ to the word~$w$ we have to rewrite it in the form~\eqref{eq:Normal-i-Form} by means of
anti-commutativity and then invert the order of $u_j$s.
Let us denote $u=\{x_{i+1}, \dots \{x_{n-1},x_n\}\dots\}$.
Then
\begin{gather*}
w = \{x_1, \dots, x_{i-1},\{x_{i},\{x_{i+1}, \dots \{x_{n-1},x_n\}\dots\}\}\dots \}
\\
\phantom{w}
= \{x_1, \dots, x_{i-1},\{x_{i}, u\}\dots \} = - \{x_1, \dots, x_{i-1},\{u, x_i\}\dots \},
\\
F_i w= -(-1)^i \{u,\{x_{i-1},\dots, \{x_1,x_i\}\dots \}\} = (-1)^i \{u,\{x_{i-1},\dots, \{x_2,\{x_i,x_1\}\}\dots\}\}.
\end{gather*}
The expression for $F_iw \in \pm V_1$ is already in the form~\eqref{eq:Normal-i-Form}.
In the same way, apply $F_1$ to $F_iw$:
\begin{gather*}
F_1F_iw= (-1)^i F_1 \{u,\{x_{i-1},\dots, \{x_i,x_1\}\dots \}\}
= \{x_i,\{x_2, \dots, \{x_{i-1}, \{u,x_1 \}\dots \}\}.
\end{gather*}
It remains to apply anti-commutativity to write the result $F_1F_iw$ in the form~\eqref{eq:Normal-i-Form}:
\begin{gather*}
F_1F_iw = -\{x_i,\{x_2, \dots x_{i-1},\{x_1,\{x_{i+1}, \dots \{x_{n-1},x_n\}\dots \}\} \}\},
\end{gather*}
$i=2,\dots, n$.
Since transpositions of the form $(1i)$ generate the entire symmetric group $S_n$, the lemma is proved.
\end{proof}

\subsection[Jacobian $\operatorname{AC}$-polynomials]{Jacobian $\boldsymbol{\operatorname{AC}}$-polynomials}\label{sec2}

In this section, we describe polylinear $\operatorname{AC}$-polynomials that have a~specif\/ic property if considered as elements of the
free $\operatorname{GP}$-algebra.

Suppose $\Psi (x_1,\dots, x_n) $ is an element of the free $\operatorname{GP}$-algebra $\operatorname{GP}(X)$, $X=\{x_1,\dots, x_n\}$ which is linear
with respect to $x_n$.
We will say that $\Psi $ is a~derivation with respect to $x_n$ if
\begin{gather*}
\Psi(x_1,\dots, x_{n-1}, yz) = y\Psi(x_1,\dots, x_{n-1}, z) + z\Psi(x_1,\dots, x_{n-1}, y)
\end{gather*}
in the free $\operatorname{GP}$-algebra $\operatorname{GP}(x_1,\dots, x_{n-1}, y,z)$.

\begin{definition}
A~polylinear $\operatorname{AC}$-polynomial $\Psi=\Psi(x_1,\dots, x_n)\in V_0$ is said to be a~{\em Jacobian} if~$\Psi $ is a~derivation
with respect to each variable~$x_i$, $i=1,\dots, n$.
A~polylinear element of $\operatorname{GP}(X)$ with the same property is called a~Jacobian $\operatorname{GP}$-polynomial.
\end{definition}

For free Lie algebra considered as a~part of the free Poisson algebra, a~similar notion was considered
in~\cite{FarkasI}.
Obviously, if $n=2$ then $C_2=\{x_1,x_2\}$ is a~Jacobian $\operatorname{AC}$-polynomial.
It was shown in~\cite{FarkasI} that there are no other Jacobian Lie polynomials (up to a~multiplicative constant).
However, there exists a~Jacobian $\operatorname{AC}$-polynomial of degree~3:
\begin{gather*}
J_3 =\{\{x_1,x_2\},x_3\} + \{\{x_2,x_3\},x_1\} + \{\{x_3,x_1\},x_2\}.
\end{gather*}
The main purpose of this section is to show that $C_2$ and $J_3$ exhaust all Jacobian $\operatorname{AC}$-polynomials.

For a~generic Poisson algebra~$A$, $a\in A$, consider the linear map $\ad a: x\mapsto \{a,x\}$, $x\in A$.
The set of all such transformations $\{\ad a\,|\, a\in A\}\subset \mathrm{End}_{\Bbbk}(A)$ generates a~Lie
subalgebra $L(A)\subset \operatorname{gl} (A) = \mathrm{End}_{\Bbbk}(A)^{(-)}$.

Given $L\in L(\operatorname{GP}(x_1,\dots, x_{n-1}))$, one may easily note that $L(x_n)\in \operatorname{GP}(x_1,\dots, x_{n})$ is a~derivation with
respect to $x_n$.
Indeed, the Leibniz identity implies that $\ad u$, $u\in \operatorname{GP}(x_1,\dots, x_{n-1})$, is a~derivation with respect
to $x_n$, and the commutator of derivations is a~derivation itself.

\begin{lemma}
\label{lem:1Der_one_var}
Let $\Psi(x_1,\dots, x_n) \in \operatorname{AC}(X)\subset \operatorname{GP}(X)$ be a~polylinear element such that $\Psi $ is a~derivation with respect
to $x_n$.
Then there exists $L\in L(\operatorname{AC}(X))$ such that $\Psi = L(x_n)$.
\end{lemma}

\begin{proof}
The algebra of multiplications $M(\operatorname{AC}(x_1,\dots, x_{n-1}))\simeq \operatorname{As}(U_n)$ contains free Lie subalgebra $\mathcal L =
\operatorname{Lie}(U_n)\subset \operatorname{As}(U_n)^{(-)}$ generated by $\ad u$ for all normal words $u\in U_n$.

As $\mathcal L$ acts on $V=\operatorname{AC}(x_1,\dots, x_{n-1}, y,z)$, we have the standard $\mathcal L$-module structure on $V\otimes
V$, given~by
\begin{gather*}
a(u\otimes v) = au\otimes v + u\otimes av,
\qquad
a\in \mathcal L,
\quad
u,v\in V.
\end{gather*}
Since $\operatorname{As}(U_n)=U(\mathcal L)$ is the universal enveloping algebra of $\mathcal L$, $V\otimes V $ is also an $U(\mathcal
L)$-module given~by
\begin{gather*}
a(u\otimes v) = \sum\limits_{(a)} a_{(1)}u\otimes a_{(2)}v,
\qquad
a\in U(\mathcal L),
\quad
u,v\in V,
\end{gather*}
where $\Delta:a\mapsto \sum\limits_{(a)} a_{(1)}\otimes a_{(2)}$ is the standard coproduct in $U(\mathcal L)$.

An $\operatorname{AC}$-polynomial $\Psi (x_1,\dots, x_n)$ may be presented as $L(x_n)$ for some $L\in \operatorname{As}(U_n)$.
By def\/inition, $\Psi $ is a~derivation with respect to $x_n$ if and only if
\begin{gather*}
L(u\otimes v) = L(u)\otimes v + u\otimes L(v) \in V\otimes V
\end{gather*}
for all $u,v\in V$.
Since $V\otimes V$ is a~faithful $U(\mathcal L)$-module, we obtain $\Delta (L) = L\otimes 1 + 1\otimes L$, thus the
Friedrichs criterion for the Lie elements in $\operatorname{As}(U_n)$ implies~$L$ to be an element of $\mathcal L$, which proves the
claim.
\end{proof}

Def\/ine a~linear map
\begin{gather*}
D(\cdot, x_n; y,z): \ V_0 \to \operatorname{GP}(x_1,\dots, x_{n-1}, y,z)
\end{gather*}
as follows: given $w=W(x_n)$, $W\in M(\operatorname{AC}(x_1,\dots, x_{n-1}))$, set
\begin{gather*}
D(w, x_n; y,z) = W(yz) - yW(z) - zW(y).
\end{gather*}

A polylinear $\operatorname{AC}$-polynomial $\Psi (x_1,\dots, x_n)\in V_0$ is a~derivation with respect to $x_n$ if and only if $D(\Psi,
x_n; y,z)=0$, or, as we have noticed above,~$W$ is primitive ($\Delta (W) = 1\otimes W + W\otimes 1$).
This property of~$W$ is homogeneous in $\operatorname{As}(U_n)$, i.e.,~$W$ splits into groups of homogeneous summands $W=W_1+\dots+W_l$,
where $W_i(x_n)$ is again a~derivation with respect to $x_n$.

\begin{remark}
\label{rem:Rem1}
Suppose $\Psi=L(x_n)\in \operatorname{AC}(X) $ is a~polynomial as in Lemma~\ref{lem:1Der_one_var}.
Then the operator~$L$ belongs to $\operatorname{Lie}(\ad u_1,\dots, \ad u_k)$ for some normal words $u_1,\dots, u_k
\in U_n$ (assume~$k$ is minimal).
For every $i=1,\dots, k$ $\Psi $ must contain a~term
\begin{gather*}
\{u_{i_1},\{u_{i_2} \dots\{u_{i},x_n\}\dots \}\},
\end{gather*}
in which the word $u_i$ appears as the last entry.
\end{remark}

Without loss of generality, assume~$L$ is polylinear with respect to $\ad u_1, \dots, \ad u_k$.
As an element of the free Lie algebra it may be uniquely written as a~linear combination of $w_{i_1,\dots, i_{k-1}} =
[\ad u_{i_1}, [\ad u_{i_2}, \dots, [\ad u_{i_{k-1}}, \ad u_i]\dots]]$.
The expansion of such a~monomial in the free associative algebra $\operatorname{As}(U_n)$ contains unique term $\ad u_{i_1}
\ad u_{i_2} \dots \ad u_{i_{k-1}} \ad u_i\in \operatorname{As}(U_n)$ ending with $\ad u_i$.
These terms for dif\/ferent $w_{i_1,\dots, i_{k-1}}$ do not cancel.

\begin{corollary}
\label{cor:ACformDer}
Let $\Psi(x_1,\dots, x_n) \in \operatorname{AC}(X)$ be a~polylinear $\operatorname{AC}$-poly\-no\-mi\-al of degree~$n$ such that~$\Psi $ is a~derivation
with respect to~$x_n$.
Suppose $\Psi = \sum\limits_{w} \alpha_w w$, $w\in X^{**}$ are normal words, and
\begin{gather*}
\max\limits_{w:\alpha_w\ne 0} \operatorname{ht}(w,x_n)\ge \max\limits_{w:\alpha_w\ne 0} \operatorname{ht}(w,x_i),
\qquad
i=1,\dots, n
\end{gather*}
$(x_n$ has maximal height in $\Psi)$.
Then
\begin{gather*}
\max\limits_{w:\alpha_w\ne 0} \operatorname{ht}(w,x_n) = n-1
\end{gather*}
and thus $\Psi $ contains a~monomial of the form
\begin{gather*}
w = \{x_{1s}, \dots \{x_{(n-1)s}, x_n\}\dots \}
\end{gather*}
for some $s\in S_{n-1}$.
\end{corollary}

\begin{proof}
Assume $k<n-1$ is the maximal height of $x_n$ in $\Psi $, i.e., $\Psi $ contains a~summand of the form $\{u_{1}, \dots
\{u_{k},x_n\}\dots \}$, $k<n-1$.
Then there exists at least one $u_i$ whose degree is greater than~1.
Remark~\ref{rem:Rem1} implies that $\Psi $ contains a~summand $\alpha_w w$, where $w = \{u_{j_1},\dots\{u_{j_k},x_n\}\dots \}$, $j_k=i$,
$\alpha_w\ne 0$.
Since $\operatorname{ht}(u_i,x_j)>1$ for some $x_j$, we have $\operatorname{ht}(w,x_j)>k$, which contradicts to the condition $\operatorname{ht}(w,x_j)\le k$.
Hence, $k=n-1$.
\end{proof}

\begin{lemma}
\label{lem:JacobianInvariant}
Suppose $\Psi=\Psi(x_1,\dots, x_n)$ is a~Jacobian $\operatorname{AC}$-polynomial.
Then~$\Psi $ is invariant with respect to the action of the group $\mathcal F$ generated by all $x_i$-flips, $i=1,\dots,
n$.
\end{lemma}

\begin{proof}
Let us f\/ix $i\in \{1,\dots, n\}$.
Without loss of generality we may assume $i=n$.
By Lemma~\ref{lem:1Der_one_var}, $\Psi = L(x_n)$, where~$L$ is a~linear operator constructed by commutators of operators
$\ad u$, $u\in U_n$.
The set of all such $\ad u$ generates an associative subalgebra $\mathcal U\subset \mathrm{End}_\Bbbk V_n$,
$\mathcal U\simeq \operatorname{As}(U_n)$.

Since $\Psi $ is polylinear (with respect to~$X$),~$L$ naturally splits into a~sum of operators presented by polylinear
(with respect to $U_n$) elements of $\mathcal U$.

Consider the linear transformation $\tau $ of $\mathcal U$ given~by
\begin{gather*}
\tau: \ W \mapsto -W^*,
\qquad
W\in \mathcal U.
\end{gather*}
The map $\tau $ acts as an identity on $\operatorname{Lie}(U_n)\subset \mathcal U^{(-)}$, it follows from the obvious observation $\tau
([W_1,W_2])=[\tau(W_1),\tau(W_2)]$ for $W_1,W_2\in \mathcal U$.

By Def\/inition~\ref{defn:x_i-flip},
\begin{gather*}
F_n(\Psi) = F_n(L(x_n)) = \tau(L)(x_n) = L(x_n) = \Psi.
\end{gather*}
As $\Psi $ is invariant with respect to all f\/lips, we have $\mathcal F(\Psi) =\{\Psi \}$.
\end{proof}

\begin{lemma}
\label{lem:AlternatingSum=Lie}
Suppose $U=\{u_1,\dots, u_m\}$ is a~set, $\operatorname{As}(U) $ is the free associative algebra generated by~$U$, $\operatorname{Lie}(U)$ is the free
Lie algebra generated by~$U$, $\operatorname{Lie}(U)\subset \operatorname{As}(U)^{(-)}$.
Consider
\begin{gather*}
A_m = \sum\limits_{s\in S_m}(-1)^s u_{1s}\cdots u_{ms}.
\end{gather*}
Then $A_m \in \operatorname{Lie}(U)$ if and only if $m=1$ or $m=2$.
\end{lemma}

\begin{proof}
For $m=1,2$ it is obvious that $A_m\in \operatorname{Lie}(U)$.

Assume $m\ge 3$ and $A_m\in \operatorname{Lie}(U)$.
Consider the homomorphism $\Phi: \operatorname{As}(U)\to \wedge (\Bbbk U)$ given by $u\mapsto u$, where $\wedge (\Bbbk U)$ is the
exterior algebra of the linear space spanned by~$U$.
Note that $\Phi(A_m) = m! u_1\cdots u_m \ne 0$ in $\wedge(\Bbbk U)$.
However, $\Phi(\operatorname{Lie}(U))\subset \wedge(\Bbbk U)^{(-)}$ is a~Lie subalgebra generated by~$U$.
It is easy to see that $\wedge(\Bbbk U)^{(-)}$ is a~3-nilpotent Lie algebra, so $\Phi(\operatorname{Lie}(U))$ does not contain elements
of degree $m\ge 3$.
\end{proof}

\begin{theorem}
Let $X=\{x_1,\dots, x_n\}$, and let $\Psi =\Psi(x_1,\dots, x_n)\in \operatorname{AC}(X)$ be a~Jacobian $\operatorname{AC}$-polynomial.
Then either $n=2$ and $\Psi = \alpha C_2$, or $n=3$ and $\Psi = \alpha J_3$, where $\alpha \in \Bbbk^*$.
\end{theorem}

\begin{proof}
By Lemma~\ref{lem:JacobianInvariant} $F\Psi = \Psi$ for every $F\in \mathcal F$.
Without loss of generality we may assume that $x_n$ has the maximal height in $\Psi $ (re-numerate variables if needed).
Corollary~\ref{cor:ACformDer} implies that $\Psi $ contains a~summand of the form $\alpha w$, where $\alpha\in\Bbbk^*$,
$w = \{x_{1s}, \dots \{x_{(n-1)s}, x_n\}\dots \}$ for some $s\in S_{n-1}$.
Without loss of generality, $\alpha =1$ and $s=\mathrm{id}$.
By Lemma~\ref{lem:FlipsGenerateS_n}, $\Psi $ contains all monomials obtained from $w $ by all permutations of variables,
i.e.,
\begin{gather*}
\Psi = \sum\limits_{s\in S_{n-1}}(-1)^s \{x_{1s}, \dots \{x_{(n-1)s}, x_n\}\dots \} + \Phi(x_1,\dots, x_n),
\end{gather*}
where the $x_n$-height of all monomials in $\Phi $ is smaller than $n-1$.
Since all summands of $\Psi $ with the same $x_n$-height form a~derivation with respect to $x_n$, the $\operatorname{AC}$-polynomial
\begin{gather*}
\Psi_1 = \sum\limits_{s\in S_{n-1}}(-1)^s \{x_{1s}, \dots \{x_{(n-1)s}, x_n\}\dots \}
\end{gather*}
must be a~derivation with respect to $x_n$.
But
\begin{gather*}
\Psi_1 = A_{n-1}(u_1,\dots, u_{n-1})(x_n),
\qquad
u_i=\ad x_i,
\end{gather*}
so $A_{n-1}(u_1,\dots, u_{n-1})\in \operatorname{Lie}(u_1,\dots, u_{n-1})$.
By Lemma~\ref{lem:AlternatingSum=Lie}, $n-1\le 2$, so $n\le 3$.
Obviously, $C_2$ and $J_3$ are the only Jacobian $\operatorname{AC}$-polynomials for $n=2$ and $n=3$, respectively.
\end{proof}

\section{Identities of generic Poisson algebras}

Let~$A$ be a~$\operatorname{GP}$-algebra, and let $f\in \operatorname{GP}(x_1,\dots, x_n)$, $f\ne 0$.
As usual, we say that~$f$ is a~polynomial identity on~$A$ if for every homomorphism $\varphi: \operatorname{GP}(x_1,\dots, x_n) \to A$
we have \mbox{$\varphi(f)=0$}.
In this case we also say that~$A$ satisf\/ies the polynomial identity~$f$.

\begin{proposition}
\label{prop:GPIdent:Jacobian}
Suppose a~$\operatorname{GP}$-algebra~$A$ satisfies a~polynomial identity.
Then there exists a~polynomial identity $\Psi $ on~$A$ which is a~Jacobian $\operatorname{GP}$-polynomial.
\end{proposition}

This statement, as well as its proof, is similar to the result by Farkas~\cite{FarkasI} on polynomial identities of
Poisson algebras.

\begin{proof}
The standard linearization procedure (see, e.g.,~\cite[Chapter~1]{ZSSS}) allows to assume that~$A$ satisf\/ies
a~polylinear polynomial identity $f\in \operatorname{GP}(X)$, $X=\{x_1,\dots, x_n\}$.

As an element of $\operatorname{GP}(X)$,~$f$ may be uniquely presented as a~linear combination of $\operatorname{GP}$-monomials $w=u_1\cdots u_k $,
$u_j\in U$, where $U\subset \operatorname{AC}(X)$ is the set of normal words.
We may assume that $u_j$ are of degree two or more (if an $\operatorname{AC}$-monomial of degree one appears, e.g., $u_j=x_i$, then one
may plug in $x_i=1$ and obtain a~polylinear polynomial identity without $x_i$).
Denote by $\operatorname{FH}_i(w)$ (the Farkas height) the degree of $u_j$ in which the variable $x_i$ occurs, and let $\operatorname{FH}_i(f)$ be the
maximal of $\operatorname{FH}_i(w)$ among all $\operatorname{GP}$-monomials~$w$ that appear in~$f$ with a~nonzero coef\/f\/icient.
Finally, set
\begin{gather*}
\operatorname{FH}(f) = \sum\limits_{i=1}^n 3^{\operatorname{FH}_i(f)}.
\end{gather*}
Observe that if~$f$ is not a~derivation in $x_i$ then the derivation dif\/ference $ D(f, x_i; x_i, x_{n+1}) $ is a~nonzero
polylinear element of $\operatorname{GP}(X\cup\{x_{n+1}\})$ which has a~smaller Farkas height.
Indeed, for a~$\operatorname{GP}$-monomial~$w$ from~$f$ we have
\begin{gather*}
\operatorname{FH}_j(D(w,x_i;x_i,x_{n+1}))\le \operatorname{FH}_j(w),
\\
\operatorname{FH}_{i}(D(w,x_i;x_i,x_{n+1}))\le \operatorname{FH}_i(w)-1,
\\
\operatorname{FH}_{n+1}(D(w,x_i;x_i,x_{n+1}))\le \operatorname{FH}_i(w)-1,
\end{gather*}
which implies
\begin{gather*}
\operatorname{FH}(w)-\operatorname{FH}(D(w,x_i;x_i,x_{n+1})) \le 3^{\operatorname{FH}_i(w)} - 2\cdot 3^{\operatorname{FH}_i(w)-1} >0.
\end{gather*}
Obviously, $D(f, x_i; x_i, x_{n+1})$ is a~polynomial identity on~$A$.

Therefore, after a~f\/inite number of steps we obtain a~nonzero polynomial identity on~$A$ which is a~Jacobian
$\operatorname{GP}$-polynomial in a~larger set of variables $\widetilde X\supseteq X$.
\end{proof}

Let us recall the notion of f\/ine grading~\cite{FarkasI}.
First, given a~set~$X$, the free anti-commutative algebra $\operatorname{AC}(X)$ carries $[X^*]$-grading such that $u\in X^{**}$ has weight~$[u]$.
Next, if $w = (u_1)\cdots (u_n)\in \operatorname{GP}(X)$, $u_i\in X^{*}$, then the weight of~$w$ is $[u_1]+\cdots + [u_n]\in \Bbbk [X^*]$.
As a~result,
\begin{gather*}
\operatorname{GP}(X) = \bigoplus\limits_{p\in \mathbb Z_+ [X^*]\setminus \{0\}} \operatorname{GP}_p(X),
\end{gather*}
where $\mathbb Z_+$ stands for the set of non-negative integers, $\operatorname{GP}_p(X)$ is the space spanned by all generic Poisson
monomials of degree~$p$.
An element $f\in \operatorname{GP}_p(X)$ is said to be {\em finely homogeneous}.

\begin{proposition}
A~Jacobian $\operatorname{GP}$-polynomial $\Psi $ can be presented as a~linear combination of products of Jacobian $\operatorname{AC}$-polynomials $($on the
appropriate set of variables$)$.
\end{proposition}

\begin{proof}
Let $X=\{x_1,x_2,\dots, x_n \}$ be a~set of variables, and let $U=\{u_1,u_2,\dots \}$ be the set of normal
nonassociative words in~$X$ (with respect to some ordering), then $\operatorname{GP}(X) = \Bbbk[U]$.
For $\Psi \in \operatorname{GP}(X)$, denote by $\supp(\Psi)$ all variables from~$X$ that appear in $\Psi $ and by
$\psupp(\Psi)$ all elements from~$U$ that appear in~$\Psi $.

Suppose $f\in \operatorname{GP}(x_1,\dots, x_n)\subseteq \operatorname{GP}(X)$ is a~Jacobian $\operatorname{GP}$-polynomial.
Without loss of generality we may assume~$f$ to be f\/inely homogeneous and $f\notin \operatorname{AC}(X)$.
Proceed by induction on $|\psupp(f)|$.

Consider a~$\operatorname{GP}$-monomial~$w$ in~$f$.
Since $f\not\in \operatorname{AC}(X)$, there exist $u_i$ and $w'\neq 1$ for which $w = u_iw'$.
Write $f = u_ig + h$, $g,h\in \operatorname{GP}(X)$, $g\ne 1$, where all $\operatorname{GP}$-monomials of~$h$ are not divisible by $u_i$ (in $\Bbbk[U]$).
Since~$f$ is polylinear, $\supp(g) \cap \supp(u_i) = \varnothing$.

Denote by $D_i$ a~map $ \operatorname{GP}(X) \to \operatorname{GP}(X\cup \{y,z\})$ def\/ined as follows:
\begin{gather*}
D_i(\Psi) =
\begin{cases}
D(\Psi, x_i; y,z), & x_i\in \supp (\Psi),
\\
\Psi, & x_i\notin \supp (\Psi).
\end{cases}
\end{gather*}
Then $D_j(f) = u_iD_j(g) + D_j(h) = 0$ if $x_j \in \supp(g)$.

Consider $\operatorname{GP}(X\cup \{y,z\})$ as a~polynomial algebra with a~set $\widetilde U$ of generators including~$U$.
Then $u_i\notin \psupp(h)$ and $u_i\notin \psupp(D_j(h))$.
Hence, $D_j(g) = 0$ and~$g$ is a~Jacobian $\operatorname{GP}$-polynomial.

Let us now f\/ix the deg-lex order on the set $[U^*]$, i.e., commutative monomials in~$U$ are f\/irst compared by their
length and then lexicographically, assuming $u_1<u_2<\cdots $.
Recall that $f=u_ig+h$, where $\psupp(h)\not\ni u_i$, and presented~$h$ as $h=gp +r$, where all $\operatorname{GP}$-monomials
of~$r$ are not divisible (in $\Bbbk[U]$) by the leading $\operatorname{GP}$-monomial $\bar g$ of~$g$.
Then $f = gq + r$, $q = u_i + p$, and $\psupp(r)\not\ni u_i$.
In particular, $\psupp(r)\subset \psupp (f)$.

By def\/inition, $D_j(f) = gD_j(q) + D_j(r) = 0$ if $x_j \in \supp(q)$.
If $D_j(q) \neq 0$ then some of the monomials in $D_j(r)$ are divisible by~$\bar g$.
Consider a~$\operatorname{GP}$-monomial~$M$ of~$r$.
Since it is not divisible by $\bar g$ there is at least one variable $u_a$ which appears in $\bar g$ and does not appear
in~$M$.
Note that if $\supp(u_b) \not\ni x_i$ then $D_i(u_b) = u_b$, and if $\supp u_b \ni x_i$ then
$D_i(u_b)$ is a~$\operatorname{GP}$-polynomial of degree two (in $\Bbbk[\widetilde U]$) in which none of variables belongs to~$U$.
Hence, $D_j(M)$ is not divisible by $u_a$ and none of the $\operatorname{GP}$-monomials of $D_j(r)$ is divisible by $\bar g$.
Therefore, $D_j(q) = 0$ and~$q$ is a~Jacobian $\operatorname{GP}$-polynomial.

Since a~product of two Jacobian $\operatorname{GP}$-polynomials is also Jacobian (with respect to the corresponding sets of variables),
$r=f-gq$ is a~Jacobian $\operatorname{GP}$-polynomial.
By induction, the statement holds for~$r$, as well as for~$g$ and~$q$.
\end{proof}

\begin{corollary}
\label{cor:2DerLeftJustified}
Let $F(t_1,\dots, t_n) \in \operatorname{GP}(t_1,t_2, \dots)$ be a~finely homogeneous Jacobian $\operatorname{GP}$-poly\-no\-mi\-al.
Then~$F$ contains a~summand $\alpha u_1\cdots u_k$, where $\alpha \in \Bbbk^*$, $u_i\in \operatorname{AC}(t_1,t_2, \dots)$ are of the
form
\begin{gather*}
\{t_{i_1},t_{i_2}\}
\qquad
\text{or}
\qquad
\{t_{i_1},\{t_{i_2},t_{i_3}\}\}.
\end{gather*}
\end{corollary}

\section{The Freiheitssatz for (generic) Poisson algebras}

The following statement is well-known in the theory of dif\/ferential f\/ields~\cite{Kolchin, Ritt}.
We will sketch a~proof below in order to make the exposition more convenient for a~reader.
Recall that the characteristic of the base f\/ield $\Bbbk $ is assumed to be zero, and that $\operatorname{Dif}_n$ denotes the variety of
commutative associative algebras with~$n$ pairwise commuting derivations.

\begin{theorem}
Every algebra from $\operatorname{Dif}_n$ which is a~field can be embedded into an algebraically closed algebra in $\operatorname{Dif}_n$.
\end{theorem}

\begin{proof}
Let~$F$ be a~dif\/ferential f\/ield of characteristic zero with a~set $\Delta=\{\partial_i \,|\, i=1,\dots, n\} $ of pairwise
commuting derivations.
Denote by $F[x;\Delta]=F*_{\operatorname{Dif}_n} \operatorname{Dif}_n(x)$ the set of all dif\/ferential polynomials in one variable~$x$ over~$F$.
Suppose $f(x)\in F[x;\Delta]\setminus F$.
Then there exists a~dif\/ferential f\/ield~$K$ which is an extension of the dif\/ferential f\/ield~$F$ such that the equation
$f(x)=0$ has a~solution in~$K$.

Indeed, dif\/ferential polynomials $F[x;\Delta]$ may be considered as ordinary polynomials in inf\/ini\-te\-ly many variables
\begin{gather*}
X=\big\{x^{(i_1,\dots, i_n)}\,\big|\, (i_1,\dots, i_n)\in \mathbb Z_+^n\big\},
\end{gather*}
where $x^{(i_1,\dots, i_n)}$ is identif\/ied with $\partial_1^{i_1}\cdots \partial_n^{i_n}(x)$.
Then the dif\/ferential ideal $I(f;\Delta)$ generated by~$f(x)$ in $F[x;\Delta]$ coincides with the ordinary ideal in~$F[X]$ generated by~$f$ and all its derivatives $\partial_1^{i_1}\cdots \partial_n^{i_n}(f)$.

Note that if $f\notin F$ then $I(f;\Delta)$ is proper: one may apply the notion of a~characteristic set (see,
e.g.,~\cite[Chapter~I.10]{Kolchin}) or simply note that the set of all derivatives of~$f$ is a~Gr\"obner basis provided that
we choose an ordering of monomials in such a~way that highest derivative (leader) is contained in the leading monomial
(e.g., rank ordering in~\cite[Chapter~I.8]{Kolchin}).
Indeed, if $uy$ is the leading monomial of~$f$ ($y\in X$ is the leader of~$f$,~$u$ is an ordered monomial in~$X$) then
$uy^{(i_1,\dots, i_n)}$ is the leading monomial of $\partial_1^{i_1}\cdots \partial_n^{i_n}(f)$.
It is easy to see that there are no compositions (we follow the terminology of Shirshov~\cite{Sh62}, see~\cite{BK2001}
for details) among~$f$ and its derivatives except for the case when $uy=y^k$, but in the latter case the only series of
compositions of intersection of~$f$ with itself is obviously trivial.

Hence, if $f\notin F$ then $I=I(f;\Delta)$ is proper, and so is its radical $\sqrt{I}$.
By the dif\/ferential prime decomposition theorem (see, e.g.,~\cite[Chapter~1]{Ritt}), $ I = p_1\cap \dots \cap p_k$, where
$p_i$ are prime dif\/ferential ideals in $F[x;\Delta]$.
In particular, $f\in p_1$, and $F[x;\Delta]/p_1$ is a~dif\/ferential domain containing a~root $x+p_1$ of~$f$.
Finally, the quotient f\/ield of that domain $Q(F[x;\Delta]/p_1)$ is the desired dif\/ferential f\/ield.

Therefore, every nontrivial equation over an arbitrary dif\/ferential f\/ield~$F$ has a~solution in an extension~$K$ of~$F$.
If~$F$ is inf\/inite then~$K$ has the same cardinality as~$F$, so the standard transf\/inite induction arguments similar to
those applied to ordinary f\/ields show that~$F$ can be embedded into a~dif\/ferential f\/ield $\bar F\in \operatorname{Dif}_n$ in which
every nontrivial dif\/ferential polynomial has a~root.
\end{proof}

\begin{corollary}[\cite{MLU10}]
The Freiheitssatz holds for the variety of Poisson algebras.
\end{corollary}

\begin{proof}
Let $A_{2n} = \Bbbk(x_1,y_1, x_2,y_2, \dots,x_n, y_n)$ be the algebra of (commutative) rational functions over $\Bbbk$,
$\partial_i= \partial_{x_i}$, $\partial'_i=\partial_{y_i}$ be ordinary partial derivatives with respect to $x_i$, $y_i$,
respectively.
As $A_{2n}\in \operatorname{Dif}_{2n}$, there exists its algebraically closed extension $\bar A_{2n} \in \operatorname{Dif}_{2n}$.
Let $\operatorname{PS}_n = A_{2n}^{(\partial)}$ be the Poisson algebra def\/ined by~\eqref{eq:PoisBracket}.
Then $\operatorname{PS}_n \subseteq \bar A_{2n}^{(\partial)}$, where the latter is a~conditionally closed Poisson algebra by Proposition~\ref{prop:CCFunctor}.

It was shown in~\cite{FarkasI} that for every nonzero Poisson polynomial $h=h(x_1,\dots, x_m)$, $m\ge 1$, there exists
a~suf\/f\/iciently large~$N$ such that $\operatorname{PS}_N$ (and thus $\bar A_{2N}^{(\partial)}$) does not satisfy the identity
$h(x_1,\dots, x_m)=0$.
Lemma~\ref{lem:CC-Frei} implies the claim.
\end{proof}

Let us twist the functor $\partial: \operatorname{Dif}_{2n}\to \operatorname{Pois}$ in order to obtain a~conditionally closed generic Poisson algebra
that does not satisfy a~f\/ixed polynomial identity.

Consider the variety $\operatorname{CDif}_{n}$ of commutative dif\/ferential algebras with pairwise commuting derivations $\partial_i$
and constants $c_i$, $i=1,\dots, n$, such that $\partial_i(c_j)=\delta_{ij}$.
Then there exists a~natural forgetful functor $\omega: \operatorname{CDif}_n\to \operatorname{Dif}_n$ erasing the information about constants.

In particular, $A_{2n}$ may be considered as an algebra from $\operatorname{CDif}_{2n}$ with derivations $\partial_i=\partial_{x_i}$,
$\partial'_i=\partial_{y_i}$, and constants $c_i=x_i$, $c_i'=y_i$, $i=1,\dots, n$.
Moreover, if $A_{2n}\subseteq A\in \operatorname{Dif}_{2n}$ then $A=B^{(\omega)}$ for an appropriate $B\in \operatorname{CDif}_{2n}$.
Hence (see Remark~\ref{rem:CCAlgebras_const}), for every $n\ge 1$ there exists a~conditionally closed algebra $\bar
B_{2n}$ in $\operatorname{CDif}_{2n}$, $\bar B_{2n}^{(\omega)}=\bar A_{2n}$.

Suppose $B\in \operatorname{CDif}_{2n}$ with derivations $\partial_i$, $\partial_i'$ and constants $c_i$, $c_i'$, $i=1,\dots, n$.
Let us consider the following functor $\tau $ from $\operatorname{CDif}_{2n}$ to the variety $\operatorname{NDif}_{2n}$ of commutative dif\/ferential
algebras with non-commuting derivations $\xi_i$, $\xi_i'$, $i=1,\dots, n$.
On the same space~$B$, def\/ine new derivations~by
\begin{gather*}
\xi_i(a) = c'_{i+1}\partial_i,
\qquad
i=1,\dots, n-1,
\\
\xi_n(a) = c'_1\partial_n,
\\
\xi'_i(a) = \partial'_i(a),
\qquad
i=1,\dots, n,
\end{gather*}
for $a\in B$.
If~$B$ is conditionally closed in $\operatorname{CDif}_{2n}$ then $B^{(\tau)}$ is conditionally closed in $\operatorname{NDif}_{2n}$.

Finally, def\/ine a~functor $\xi $ from $\operatorname{NDif}_{2n}$ to the variety $\operatorname{GP}$ of generic Poisson algebras by means of
\begin{gather*}
\{a,b\} = \sum\limits_{i\ge 1} \xi_i(a)\xi'_i(b) - \xi_i(b)\xi'_i(a).
\end{gather*}
Denote by $\operatorname{GPS}_n$ the $\operatorname{GP}$-algebra $\big(A_{2n}^{(\tau)} \big)^{(\xi)}$.

\begin{proposition}
\label{prop:NonPIGPS}
For every $n\ge 1$ there exists $N\ge 1$ such that the $\operatorname{GP}$-algebra $\operatorname{GPS}_{m}$ does not satisfy a~polynomial identity of
degree~$n$ for all $m\ge N$.
\end{proposition}

\begin{proof}
Suppose $f \in \operatorname{GP}(t_1,t_2,\dots)$ is a~$\operatorname{GP}$-polynomial of degree~$n$ which is an identity on $\operatorname{GPS}_{m}$.
By Proposition~\ref{prop:GPIdent:Jacobian} there also exists a~polylinear identity~$g$ on $\operatorname{GPS}_{m}$ which is a~Jacobian
$\operatorname{GP}$-polynomial.

Let us split~$g$ into f\/inely homogeneous components:
\begin{gather*}
g = g_1 + \cdots + g_k,
\end{gather*}
each $g_i$ is a~Jacobian $\operatorname{GP}$-polynomial (but not an identity on $\operatorname{GPS}_{m}$).

According to Corollary~\ref{cor:2DerLeftJustified}, $g_1$ contains a~summand $\alpha u_1\cdots u_l$, $\alpha \in
\Bbbk^*$,
\begin{gather*}
u_i = \{t_{i_1}, \dots \{t_{i_{m_i}}, t_{i_{m_i+1}} \}\dots \},
\qquad
m_i=1,2.
\end{gather*}
Assume~$m$ is large enough (e.g., $m>2 l$), and evaluate the variables in such a~way that
\begin{gather*}
t_{i_{m_i+1}} = y_{k_i},
\quad
t_{i_{m_i}}=x_{k_i},
\quad
t_{i_{m_i-1}}=x_{k_i+1},\quad \dots,\quad  t_{i_1}=x_{k_i+m_i-1},
\\
k_{i+1}\ge k_i+m_i,
\qquad
k_l+m_l<m.
\end{gather*}
Then the only summand in $g_1(t_1,\dots, t_n) $ is nonzero, namely, the summand mentioned in
Corollary~\ref{cor:2DerLeftJustified}: it turns into $\alpha y_{k_1+m_1}\cdots y_{k_l+m_l} \ne 0$.
Other $g_i's$ turn into zero.

Hence,~$g$ cannot be a~polynomial identity on $\operatorname{GPS}_{m}$.
\end{proof}

\begin{theorem}
The Freiheitssatz holds for the variety of generic Poisson algebras.
\end{theorem}

\begin{proof}
Given $N\ge 1$, $G_N = \big(\bar B_{2N}^{(\tau)}\big)^{(\xi)}$ is a~conditionally closed algebra in $\operatorname{GP}$~by
Proposition~\ref{prop:CCFunctor}, and $\operatorname{GPS}_N\subseteq G_N$.
The claim now follows from Proposition~\ref{prop:NonPIGPS} and Lemma~\ref{lem:CC-Frei}.
\end{proof}

\subsection*{Acknowledgements}
P.~Kolesnikov was partially supported by FAPESP (grant 2012/04704-0) and RFBR (12-01-33031), L.~Makar-Limanov
acknowledges the support of FAPESP (grant 2011/52030-5), and I.~Shestakov was supported by FAPESP (grant 2010/50347-9)
and CNPq (grant 3305344/2009-9).
L.~Makar-Limanov and I.~Shestakov gratefully acknowledge Max-Planck-Institut f\"ur Ma\-the\-matics, L.~Makar-Limanov is also grateful to the University of Michigan and the Weizmann Institute of Science.
The authors are grateful to Alexei Ovchinnikov for valuable discussions and to the referees for important remarks that
helped to improve the exposition.

\pdfbookmark[1]{References}{ref}
\LastPageEnding

\end{document}